\documentclass[12pt]{article}
\usepackage{amsmath, amssymb, amsthm, enumerate}
\newtheorem{theorem}{Theorem}[section]
\newtheorem{lemma}[theorem]{Lemma}
\newtheorem{proposition}[theorem]{Proposition}
\newtheorem{remark}[theorem]{Remark}

\newtheorem{cor}[theorem]{Corollary}
\newtheorem{example}[theorem]{Example}
\newtheorem{defin}[theorem]{Definition}

\def\rr{{\mathbb R}}

\def\cc{{\mathbb C}}
\def\nn{{\bf N}}

\def\nn{{\mathbb N}}
\def\qq{{\mathbb Q}}

\def\al{\alpha}
\def\be{\beta}

\def\De{\Delta}

\def\ph{\phi}

\def\stb{,\ldots ,}

\def\msk{\medskip}

\def\noi{\noindent}

\def\sumin{\sum_{i=1}^n}

\def\sumj0m{\sum_{j=0}^m}

\def\sumi0n{\sum_{i=0}^n}
\def\sumina{\sum_{i=1}^n a_i}

\def\suminf{\sum_{i=1}^n a_if(\al_ix+\be_iy)}
\def\suminftilde{\sum_{i=1}^n a_i\tilde{f}(\al_ix+\be_iy)}

\begin{tiny}\end{tiny}

\def\proof{\noi {\bf Proof.} }

\def\x1n{x_1 \stb x_n}
\def\y1n{y_1 \stb y_n}

\def\diag{{\rm diag}\, }

\begin{document}
\title{On spectral analysis in varieties containing the solutions of inhomogeneous linear functional equations}
\author{Gergely Kiss \\ {gergely.kiss@uni.lu} \and Csaba Vincze \\ {csvincze@science.unideb.hu}}

\footnotetext[2]{{\bf MR subject classification:} primary 43A45, 43A70,
secondary 13F20}
\footnotetext[3]{G. Kiss is supported by the Internal Research Project R-STR-1041-00-Z of the University of Luxembourgh and by the Hungarian National Foundation for Scientific Research,
Grant No. K104178. Cs. Vincze is supported by the University of Debrecen's internal research project RH/885/2013.}

\maketitle

\abstract{The aim of the paper is to investigate the solutions of special inhomogeneous linear functional equations by using spectral analysis in a translation invariant closed linear subspace of additive/multiadditive functions containing the restrictions of the solutions to finitely generated fields. The application of spectral analysis in some related varieties is a new and important trend in the theory of functional equations; especially they have successful applications in case of homogeneous linear functional equations. The foundation of the theory can be found in M. Laczkovich and G. Kiss \cite{KL}, see also G. Kiss and A. Varga \cite{KV}. We are going to adopt the main theoretical tools to solve some inhomogeneous problems due to T. Szostok \cite{KKSZ08}, see also \cite{KKSZ} and \cite{KKSZW}. They are motivated by quadrature rules of approximate integration.  

\section{Introduction and preliminaries}\label{s1}

Let $\cc$ denote the field of complex numbers. We are going to investigate functional equation
\begin{equation}\label{e1}
F(y)-F(x)=(y-x)\suminf,
\end{equation}
where $x,y\in \cc$ and $f,F:\cc\to \cc$ are unknown functions. Equation \eqref{e1} is motivated by quadrature rules of approximate integration. The problem is due to T. Szostok \cite{KKSZ08}, see also \cite{KKSZ} and \cite{KKSZW}. To formulate the basic preliminary results and facts we need the notion of generalized polynomials. 
Let $(G, *)$ be an Abelian group;  $\cc ^G$ denotes the set of complex valued functions defined on $G$. A function $f:G\to \cc$ is a {\it
generalized polynomial}, if there is a non-negative integer $p$ such that
\begin{equation}\label{e4}
\De _{g_1}\ldots \De _{g_{p+1}}f=0
\end{equation}
for any $g_1 \stb g_{p+1} \in G$. Here $\De _g$ is the difference operator
defined by $\De _g f(x)=f(g*x)-f(x) \ (x\in G)$, where $f\in \cc ^G$
and $g\in G$. The smallest $p$ for which \eqref{e4} holds for any $g_1 \stb g_{p+1} \in G$ is the {\it degree} of the
generalized polynomial $f$. A function $F:G^p \to \cc$ is {\it $p$-additive}, if it is additive
in each of its variables. A function $f\in \cc ^G$ is called a
{\it generalized monomial of degree $p$}, if there is a
symmetric $p$-additive function
$F$ such that $f(x)=F(x\stb x)$ for any $x\in G$. It is known that any generalized polynomial function can be written as the sum of generalized monomials \cite{SZ3}. By a general result of M. Sablik \cite{S} any solution of \eqref{e1} is a
generalized polynomial under some mild conditions for the parameters in the functional equation. For the proof of the following result see also Lemma 2 in \cite{KKSZ}.

\begin{lemma}\label{l01} 
Let $n\in \nn$ be a given natural number. Suppose that
\begin{enumerate}
\item $\al_1\stb \al_n, \be_1 \stb \be_n\in \rr$ or $\ \cc$,
\item $\al_i+\be_i\ne 0$,
\item \begin{equation}\label{efelt1}\left|
\begin{array}{ll}
\al_i & \be_i \\
\al_j & \be_j
\end{array} \right| \ne 0, ~~i\ne j, ~ i, j\in \{1\stb n\}.\end{equation}
\end{enumerate}
If the functions $F, f_1\stb f_n: \rr\to \rr$ or $\cc\to \cc$ satisfy functional equation
\begin{equation}\label{e5}
F(y)-F(x)=(y-x)\sum_{i=1}^n f_i(\al_i x +\be_i y),~~x,y\in \rr \textrm{ or } \ \cc,
\end{equation}
then $f_1 \stb f_n$ are generalized polynomial functions of degree at most $2n-1$.
\end{lemma}

The first and the third conditions also appear in L. Sz\'ekelyhidi's result \cite{SZ3}. It states a similar conclusion under the special choice $F=0$, i.e. if 
$$\sum_{i=1}^n f_i(\al_ix +\be_iy)=0$$
and the real (or complex) parameters satisfy condition \eqref{efelt1} in Lemma \ref{l01}, then $f_i$'s are generalized polynomials of degree at most $n-2$. Using rational homogenity we can prove the following statement too (see Lemma 3 in \cite{KKSZ}).
\begin{lemma}\label{l1}
Suppose that the functions $f,F:\rr\to \rr$ or $\ \cc\to \cc$ satisfy 
equation \eqref{e1}. If the function $f$ is a generalized polynomial, i.e. it is of the form 
\begin{equation}\label{e7}
f(x)=\sum_{j=1}^p d_j(x),
\end{equation}
where $d_j$ is a monomial of degree $j\in \{1\stb p\}$ 
and $d_0$ is constant then the function $d_j$ also satisfies \eqref{e1} with some $j$-additive $F_j$ for any $j\in \{1\stb p\}$.
\end{lemma}

In what follows we suppose that conditions of Lemma \ref{l01} are satisfied, i.e. any solution of \eqref{e1} is a generalized polynomial of degree at most $2n-1$. In the sense of Lemma \ref{l1} we may also assume that $f$ is a monomial of degree at most $2n-1$ without loss of generality. Using the notation 
$$h(x)=\sumina f(\al_i x)$$
let us choose $y=0$ in equation \eqref{e1}: $\displaystyle{F(0)-F(x)=-xh(x)}$.
By substituting $x=0$ we have that  $\displaystyle{F(y)-F(0)=y\sumina f(\be_i y)}$. Therefore 
$$h(x)=\sumina f(\al_i x)=\sumina f(\be_i x)\ \ \Rightarrow\ \ F(y)-F(x)=yh(y)-xh(x)
$$ 
and, consequently, 
\begin{equation}
\label{hx}
(y-x)\sum_{i=1}^n f_i(\al_i x +\be_i y)=yh(y)-xh(x),~~x,y\in \rr \textrm{ or } \ \cc.
\end{equation}
The following theorem shows that $h(x)$ must be of a special form provided that $f$ is a monomial of degree $p$. 
\begin{theorem}
\label{KKT}
If the functions $f,F:\rr\to \rr$ or $\ \cc\to \cc$ satisfy 
equation \eqref{e1} and $f$ is a monomial of degree $p$ then $h(x)=c\cdot x^p$ for some constant $c\in \rr $ or $\cc$.
\end{theorem}

\proof The proof can be found in \cite{KKSZ} (Theorem 1).
\hfill $\square$
\medskip

In case of $p=1$ (additive solutions) Theorem \ref{KKT} allows us to simplify \eqref{hx} as 

%





\begin{equation}\label{e2}
\suminf=c\cdot (x+y) \qquad (x,y\in \cc ).
\end{equation}
Note that for any additive solution $f$ the system of equations
$$c\cdot x=\sumina f(\al_i x)\ \ \textrm{and}\ \ c\cdot y=\sumina f(\be_i y)$$
is equivalent to 
$$c\cdot (x+y)= \sumina f(\al_i x+\be_i y).$$ 

\begin{remark}\label{degreek} {\rm In case of a monomial solution $f$ of higher degree $p>1$ 
$$c\cdot x^p=\sumina f(\al_i x)\ \ \textrm{and}\ \ c\cdot y^p=\sumina f(\be_i y)$$
because of Theorem \ref{KKT} and the alternate substitutions $x=0$ and $y=0$. Equation \eqref{hx} reduces to}
\begin{equation}\label{ek}
\suminf=c\cdot \sum_{l=0}^p x^l y^{p-l} \qquad (x,y\in \cc ).
\end{equation}
\end{remark}

In what follows we are going to show that spectral analysis can be applied in translation invariant closed linear subspaces (varieties) of additive functions on some finitely generated fields containing the restrictions of the solutions of functional equation \eqref{e2}. Note that the translation invariance is taken with respect to the multiplicative group structure. We will use spectral analysis in some related varieties of equation \eqref{ek} too. The theory has successful applications in case of homogeneous linear functional equations. The case of $c=0$ has been investigated in \cite{KV} and \cite{KL}. 
In this paper we use the same ideas to generalize the results for $c\neq 0$ (inhomogeneous case). The spectral analysis provides the existence of nonzero exponential functions in the varieties of the restricted solutions to finitely generated fields. It is a necessary condition for the existence of nonzero solutions of the functional equation. The sufficiency also follows in some special cases which means the complete description of the solution space. In general we also need the application of spectral synthesis in the varieties to give the description of the solution space on finitely generated fields. Since its central objects are the so-called differential operators having a more complicated behavior relative to the automorphism solution, the application of the spectral synthesis is presented in a forthcoming paper \cite{KV1} as the continuation of the recent one.


\section{Varieties generated by non-trivial solutions of linear functional equations}

Let $G$ be an Abelian group. By a {\it variety} we mean a translation invariant closed linear
subspace of $\cc^G$.

\subsection{Varieties of additive solutions}

Consider the linear functional equation
\begin{equation}\label{e1o1}
\sum_{i=1}^n a_if(\al_ix+\be_iy)=c\cdot (x+y) \qquad (x,y\in \cc )
\end{equation}
and let a finitely generated subfield $K\subset \cc$ containing the
parameters $\al_i ,\be_i$ $(i=1\stb n)$ be fixed. If $V_1$ is the set of additive functions on $K$ then it is a closed linear subspace in $\cc ^{K}$ as the following lemma shows; for the proof see \cite{KL}. 

\begin{lemma}\label{lv1s0}
$V_1$ is a closed linear subspace of $\ \cc ^{K}$.
\end{lemma}

Now we are going to determine the smallest closed linear subspace in $V_1$ containing the additive solutions of \eqref{e2} on $K$ under $c\neq 0$. The key step is to consider $c$ as a free parameter running through the domain $\cc$. 

\begin{defin}
Let $S_1$ be the subset of $\ V_1$, where $\tilde{f}\in S_1$ if and only if there exists $\tilde{c}\in \cc$ such that  
\begin{equation}\label{e2tilde}
\suminftilde=\tilde{c}\cdot (x+y) \qquad (x,y\in K ).
\end{equation}
\end{defin}

\begin{lemma}\label{lv1s1}
$S_1$ is a closed linear subspace of $\ V_1$.
\end{lemma}

\proof It is clear that $S_1$ is a linear subspace in $V_1$. Let $\tilde{g}\colon K \to \mathbb{C}$ be a function in the closure of $S_1$. In the sense of Lemma \ref{lv1s0}, the function $\tilde{g}$ is additive. We are going to show that 
$$\sumina \tilde{g}(\al_ix)=\tilde{c}\cdot x\ \ \textrm{and}\ \ \sumina \tilde{g}(\be_i y)=\tilde{c}\cdot y$$
for some $\tilde{c}\in \cc$. Since $\tilde{g}$ is in the closure of $S_1$, for any 
$\varepsilon>0$ there exists $\tilde{f}_{\varepsilon} \in S_1$ such that
$$|\tilde{f}_{\varepsilon}(\al_i x)-\tilde{g}(\al_i x)|< \varepsilon\ \ \textrm{and}\ \ |\tilde{f}_{\varepsilon}(\be_i y)-\tilde{g}(\be_i y)|< \varepsilon$$ 
for any element $x, y$ of fixed finite subsets $X, Y\subset K$ and  $\al_i, \be_i\in K$. On the other hand  
\begin{equation}\label{ecxe} \sum_{i=1}^n a_i \tilde{f}_{\varepsilon}(\al_ix)=\tilde{c}_{\varepsilon}\cdot
x \ \ \textrm{and}\ \ \sum_{i=1}^n a_i \tilde{f}_{\varepsilon}(\be_i y)=\tilde{c}_{\varepsilon}\cdot
y 
\end{equation} 
for some $\tilde{c}_{\varepsilon} \in \cc$, respectively. First we show that there exists a finite limit
$\tilde{c}:=\lim_{\varepsilon\to 0} \tilde{c}_{\varepsilon}$. If $x=1$ then the $\De$-inequality says that
\begin{equation}
\begin{split}
&|\tilde{c}_{\varepsilon_1}-\tilde{c}_{\varepsilon_2}|=|\sumina \tilde{f}_{\varepsilon_1}(\al_i)-\sumina \tilde{f}_{\varepsilon_2}(\al_i)| \le\\
&|\sumina \tilde{f}_{\varepsilon_1}(\al_i)-\sumina \tilde{g}(\al_i)|+|\sumina \tilde{f}_{\varepsilon_2}(\al_i)-\sumina \tilde{g}(\al_i)|\le\\
&\sumin |a_i|\cdot | \tilde{f}_{\varepsilon_1}(\al_i)-\tilde{g}(\al_i)|+ \sumin |a_i|\cdot | \tilde{f}_{\varepsilon_2}(\al_i)-\tilde{g}(\al_i)|\le \lambda\cdot (\varepsilon_1 +\varepsilon_2),
\end{split}
\end{equation}
where $\lambda= \sumin|a_i|$. By Cauchy criterion, this implies that there exists
$\lim_{\varepsilon\to 0} \tilde{c}_{\varepsilon}=\tilde{c}<\infty$. On the other hand  
\begin{equation}
\begin{split}
&|\sumina \tilde{g}(\al_ix)-\tilde{c}\cdot x|\le \\
&|\sumina \tilde{g}(\al_ix)-\sumina \tilde{f}_{\varepsilon}(\al_i x)|+
|\sumina \tilde{f}_{\varepsilon}(\al_ix)-\tilde{c}\cdot x|\le \\
&|\sumina \tilde{g}(\al_ix)-\sumina \tilde{f}_{\varepsilon}(\al_i x)|+
|\sumina \tilde{f}_{\varepsilon}(\al_ix)-\tilde{c}_{\varepsilon} \cdot x |+\\
&|\tilde{c}_{\varepsilon} \cdot x-\tilde{c}\cdot x| \le \lambda \cdot \varepsilon+|\tilde{c}_{\varepsilon}-\tilde{c}|\cdot |x|
\end{split}
\end{equation}
and, in a similar way, 
\begin{equation}
|\sumina \tilde{g}(\be_i y)-\tilde{c}\cdot y|\le \lambda \cdot \varepsilon+|\tilde{c}_{\varepsilon}-\tilde{c}|\cdot |y|.
\end{equation}
Taking the limit $\varepsilon\to 0$ we have 
$$\sumina \tilde{g}(\al_ix)=\tilde{c}\cdot x\ \ \textrm{and}\ \ \sumina \tilde{g}(\be_i y)=\tilde{c}\cdot y$$
for any element $x, y$ of finite subsets $X, Y\subset K$ and  $\al_i, \be_i\in K$. Therefore they hold for any $x, y\in K$. Since by Lemma \ref{lv1s0} the function $\tilde{g}$ is additive it follows that $\tilde{g}\in S_1$. \hfill $\square$
\msk

\begin{remark} {\rm The analogue results for $c=0$ can be found in \cite{KL}.}
\end{remark}

Let $K^*= \{x\in K: x\neq 0\}$ be the Abelian group
with respect to the multiplication in $K$. We also put $V_1^* =\{ f|_{K^*} : f\in V_1 \}$ and
$$S_1^* =\{ \tilde{f}|_{K^*} : \tilde{f}\in S_1 \} .$$

\begin{lemma} \label{ls1}
$V_1^{*}$ and $S_1^*$ are varieties in $\cc^{K^*}$.
\end{lemma}

\proof It is easy to see that $V_1^{*}$ is a variety. Since $S_1$ is
a closed linear subspace of $V_1$ we have that $S_1^*$ is a closed linear
subspace of $V_1^*$. \hfill $\square$
\msk

 Recall that the translation invariance is taken with respect to the multiplicative group structure, i.e. if $\tilde{f}\in S^*_1$, then the map $\tau_a \tilde{f}\colon x\in K^* \mapsto \tilde{f}(a x)$ also belongs to $S_1^*$ for every  $a \in  K^*$.

\begin{defin} $S_1^0$ is the subspace of $S_1$ belonging to the homogeneous case $\tilde{c}=0$.
\end{defin}

\subsection{Varieties generated by higher order monomial solutions}

Consider the linear functional equation

\begin{equation}\label{eke}
\suminf=c_p \cdot \sum_{l=0}^p x^l y^{p-l} \qquad (x,y\in \cc )
\end{equation}
and let a finitely generated subfield $K\subset \cc$ containing the
parameters $\al_i ,\be_i$ $(i=1\stb n)$ be fixed. If $V_p$ is the set of $p$-additive functions on $K$ then it is a closed linear subspace in $\cc^{G}$, where $G=K\times \ldots \times K$ is the Cartesian product ($p$-times) of $K$ with itself; for the proof see \cite{KL}. 

\begin{lemma}\label{lvps0}
$V_p$ is a closed linear subspace of $\ \cc ^{G}$.
\end{lemma}

Suppose that 
$$f(x):=F_p(x, \ldots, x)$$
is a monomial solution of degree $p>1$ of equation 
\begin{equation}\label{klon}
\suminf=c_p\cdot \sum_{l=0}^p x^l y^{p-l} \qquad (x,y\in K ),
\end{equation}
where $F_p$ is a symmetric, $p$-additive function; cf. equation \eqref{ek}. By the "binomial theorem" equation \eqref{klon} implies that 
\begin{equation}
\label{binomial}\sum_{i=1}^na_i \sum_{l=0}^p\binom{p}{l}F_{p, l} (\al_i x, \be_iy)=c_p\cdot \sum_{l=0}^p x^l y^{p-l},
\end{equation}
where the script $l$ shows the number of appearences of the term $\al_i x$ among the arguments of $F_p$:
$$F_{p, l} (\al_i x, \be_iy):=F_k(\ \underbrace{\al_i x, \ldots, \al_ix}_{l \ -\ \textrm{times}}, \be_iy, \ldots, \be_iy).$$
Since additivity implies rational homogenity, rational substitutions allow us to compare the sides of equation \eqref{binomial} member by member as classical polynomials:
\begin{equation}\label{multi}
\begin{split}
&\sum_{i=1}^n a_i F_{p, p} (\al_i x, \be_iy)=\sum_{i=1}^n a_i F_p(\al_ix, \ldots, \al_i x)=c_p\cdot x^p,\\
&\sum_{i=1}^n a_i F_{p, 0} (\al_i x, \be_iy)=\sum_{i=1}^n a_i F_p(\be_iy, \ldots, \be_i y)=c_p\cdot y^p,\\
&\sum_{i=1}^n a_i \binom{p}{l} F_{p, l} (\al_i x, \be_iy)=c_p\cdot x^l y^{p-l}\ \ (l=1, \ldots, p-1).
\end{split}
\end{equation}
It is well-known that the diagonals determine the symmetric multiadditive functions. Therefore we can write conditions (\ref{multi}) into the multivariable form
{\footnotesize{\begin{equation}\label{multimulti}
\begin{split}
&\sum_{i=1}^n a_i F_p(\al_ix_1, \ldots, \al_i x_p)=c_p\cdot x_1\cdot \ldots \cdot x_p,\\
&\sum_{i=1}^n a_i F_p(\be_i y_1, \ldots, \be_i y_p)=c_p\cdot y_1 \cdot \ldots \cdot y_p,\\
&\sum_{i=1}^n a_i \binom{p}{l}F_{p} (\al_i x_1, \ldots, \al_i x_l, \be_i y_1, \ldots, \be_i y_{p-l})=c_p\cdot x_1\cdot \ldots \cdot x_l\cdot y_1 \cdot \ldots \cdot y_{p-l},
\end{split}
\end{equation}}}

\noindent
where $l=1, \ldots, p-1$; for the similar trick see \cite{KVV}. The multivariable form of the conditions will have an important role to clarify the consequences of \eqref{multi} for the properties of the translated mapping 
$$\tau_{\vec{z}}F_p (w_1, \ldots, w_p):=F_p(z_1 w_1, \ldots, z_p w_p).$$
Note that the symmetry condition is taking to fail by the translation invariance with respect to the multiplicative group structure because it is not a symmetric mapping in general. 

\begin{defin} For any $p$-additive function $F_p$ let us define $F_p^{\sigma}$ as
$$F_p^{\sigma}(w_1, \ldots, w_p):=F_p(w_{\sigma(1)}, \ldots, w_{\sigma(p)}),$$
where $\sigma$ is a permutation of the elements $1, \ldots, p$.
\end{defin}

\begin{lemma}\label{trans} Suppose that $f(x):=F_p(x, \ldots, x)$
is a monomial solution of degree $p>1$ of equation \eqref{klon}, where $F_p$ is a symmetric, $p$-additive function. For any permutation $\sigma$ of the elements $1, \ldots, p$ the mapping $\left(\tau_{\vec{z}}F_p\right)^{\sigma}$ satisfies condition \eqref{multimulti} with $\tilde{c}_p=c_p\cdot z_1\cdot \ldots \cdot z_p$ and the diagonalization 
$$\tilde{f}_p(x):=\tau_{\vec{z}}F_p (x, \ldots, x)$$
is a solution of functional equation 
\begin{equation}\label{klon1}
\sum_{i=1}^n a_i f(\al_ix+\be_iy)=\tilde{c}_p\cdot \sum_{l=0}^p x^l y^{p-l} \qquad (x,y\in K ).
\end{equation}
\end{lemma}
The proof is a straightforward calculation. Now we are in the position to define the variety generated by the higher order solutions.

\begin{defin} Let $S_p$ be the subset of $V_p$, where $\tilde{F}_p\in S_p$ if and only if there exists $\tilde{c}_p\in \cc$ such that  
{\footnotesize{\begin{equation}\label{multimultitilde}
\begin{split}
&\sum_{i=1}^n a_i\tilde{F}_p^{\sigma}(\al_ix_1, \ldots, \al_i x_p)=\tilde{c}_p\cdot x_1\cdot \ldots \cdot x_p,\\
&\sum_{i=1}^n a_i\tilde{F}_p^{\sigma}(\be_i y_1, \ldots, \be_i y_p)=\tilde{c}_p\cdot y_1 \cdot \ldots \cdot y_p,\\
&\sum_{i=1}^n a_i \binom{p}{l} \tilde{F}_{p}^{\sigma} (\al_i x_1, \ldots, \al_i x_l, \be_i y_1, \ldots, \be_i y_{p-l})=\tilde{c}_p\cdot x_1\cdot \ldots \cdot x_l\cdot y_1 \cdot \ldots \cdot y_{p-l}
\end{split}
\end{equation}}}

\noindent
$(l=1, \ldots, p-1)$ for any permutation $\sigma$ of the elements $1, \ldots, p$. 
\end{defin}

The condition \emph{for any permutation of the variables} substitutes the symmetry condition in the following sense: $S_p$ becomes closed under the usual symmetrization process
$$\left(\textrm{Sym}\ \tilde{F}_p\right) (w_1, \ldots, w_n):=\frac{1}{p!}\sum_{\sigma} \tilde{F}_p^{\sigma}(w_1, \ldots, w_p),$$
where $\sigma$ runs through the permutations of $1, \ldots, p$, i.e. the  diagonals of the elements in $S_p$ are the solutions of functional equations of type \eqref{klon1}. The problem disappears in case of $p=1$ because there are neither mixed blocks nor positions for variables to set free: $\displaystyle{\sum_{i=1}^n a_i\tilde{F}_1(\al_ix)=\tilde{c}_1\cdot x}$ and $\displaystyle{\sum_{i=1}^n a_i\tilde{F}_1(\be_i y)=\tilde{c}_1\cdot y}$. By adopting the arguments of subsection 2.1 to the multivariate setting on a finitely generated field we can prove the analogue statements of Lemma \ref{lv1s1} and Lemma \ref{ls1}.

\begin{lemma}\label{lvpsp}
$S_p$ is a closed linear subspace of $\ V_p$.
\end{lemma}

The sets $V_p^*$ and $S_p^*$ consist of the restrictions of the elements in $V_p$ and $S_p$ to the multiplicative group
$$G^*:=\underbrace{K^*\times \ldots \times K^*}_{\textrm{$p$-times}},$$
respectively. 

\begin{lemma} \label{lsp}
$V_p^{*}$ and $S_p^*$ are varieties in $\cc^{G^*}$.
\end{lemma}

\begin{remark} \rm{Lemma \ref{trans} shows that $S_p^*$ is the variety in $\cc^{G^*}$ generated by the restrictions of symmetric $p$-additive functions to $G^*$ provided that the diagonalizations are the solutions of functional equation \eqref{klon1} for some $\tilde{c}_p\in \cc$. Note that the translation invariance is taken with respect to the multiplicative group structure of $G^*$.} 
\end{remark}

\begin{defin} $S_p^0$ is the subspace of $S_p$ belonging to the homogeneous case $\tilde{c}_p=0$.
\end{defin}

\subsection{Spectral analysis in the varieties on a discrete Abelian group with torsion free rank less than continuum}

Let $(G, *)$ be an Abelian group. A function $m\colon G\to \cc$ is called \emph{exponential} if it is multiplicative:
$$m(x*y)=m(x)m(y)$$
for any $x, y\in G$. If a variety contains an exponential function  
then we say that {\it
spectral analysis holds in this variety}. If spectral analysis holds
in each variety on $G$, then {\it spectral analysis
holds on $G$.} Let $r_0 (G)$ be the torsion free rank of $G$. The following
theorem is the main result of \cite{LSZG}.
\begin{theorem}\label{tspanal}
Spectral analysis holds on a discrete Abelian group $G$ if and only
if $r_0(G)<2^{\omega}$.
\end{theorem}

\section{Applications of spectral analysis I}

The following proposition \cite[Theorem 14.5.1, p. 358]{K}) will be frequently used.

\begin{proposition}\label{p1} Let $K\subset\mathbb{C}$  be a finitely generated field and $\phi\colon K\to \mathbb{C}$
be an injective homomorphism. Then there exists an automorphism
$\psi$ of $\ \mathbb{C}$ such that $\psi|_K=\phi.$
\end{proposition}

As a consequence of Proposition \ref{p1} we can also show the following lemma
which can be found in \cite{KV}; for the definition of $K^*$ and $S_1^*$ see subsection 2.1.

\begin{lemma}\label{l30}
If $m\in S_1^*$ is an exponential on $K^*$, then there is an extension of $m$ to $\cc$ as an automorphism.
\end{lemma}

\subsection{Additive solutions of linear functional equati\-ons}

Using Theorem \ref{tspanal} spectral analysis holds in the variety of $S_1^*$ containing the restrictions of additive solutions of equation \eqref{e2} as $c$ is running through $\cc$. Therefore we can formulate the following result for the existence of a nontrivial additive solution. 

\begin{theorem}\label{t1o1}
The existence of a nonzero additive solution of \eqref{e2} implies that there exist a finitely generated subfield $K\subset \cc$ containing $\al_i$ and $\be_i$ $(i=1, \ldots, n)$ and an automorphism $\phi\colon\mathbb{C}\to\mathbb{C}$ as the extension of an exponential element in $S_1^*$ such that 
\begin{equation}\label{eA1}
\sum_{i=1}^n a_i\phi(\al_ix+\be_iy)=\tilde{c}\cdot (x+y)\ \ (x, y\in K)
\end{equation}
for some $\tilde{c}\in \cc$. Especially,
\begin{equation}\label{eA}
\sum_{i=1}^n a_i\phi(\al_i)=\sum_{i=1}^n a_i\phi (\be_i)=\tilde{c}.
\end{equation}
If $\tilde{c}=0$ then 
\begin{equation}\label{eA0}
\sum_{i=1}^n a_i\phi(\al_ix+\be_iy)=0 \ \ (x, y\in \cc),
\end{equation}
i.e. $\phi$ is the solution of the homogeneous equation on $\cc$. If $\tilde{c}\neq 0$ then $\phi(x)=x$ $(x\in K)$ and
\begin{equation}
\label{simplified}
\sum_{i=1}^n a_i\al_i=\sum_{i=1}^n a_i\be_i=\tilde{c}\neq 0.
\end{equation}
Conversely, if \eqref{simplified} holds then $f:=(c/\tilde{c})\cdot x$ is a nonzero particular additive solution of \eqref{e2} on $\cc$.
\end{theorem}

\proof
Suppose that $f$ is a nonzero additive solution of \eqref{e2}, i.e. $f(e)\neq 0$ for some $e\in\mathbb{C}$. Let $K=\mathbb{Q} (\al_1,\ldots , \al_n, \be_1, \ldots ,\be_n, e)$ be the extension of $\mathbb{Q}$ by the
complex numbers $\al_i, \be_i$ and $e$ ($i=1,\ldots ,n $).
In the sense of Lemma \ref{ls1}, $S_1^*$ is a variety in $\cc^{K^*}$. We have that $S_1^*\ne \{ 0\}$ because $f |_{K^*} \in S_1^*$ and $f(e)\neq 0$. Since $K^*$ is countable we find, by Theorem \ref{tspanal}, that $S_1^*$ contains an exponential element $\phi$, i.e. 
$$\sum_{i=1}^n a_i \phi(\al_i x+\be_i y)=\tilde{c}\cdot (x+y)\ \ (x, y\in K)$$
for some $\tilde{c}\in \cc$. By Lemma \ref{l30}, $\phi$ can be extended to an
automorphism of $\mathbb{C}$. Especially   
\begin{equation}
\label{key1}
\tilde{c}\cdot x=\sum_{i=1}^n a_i \phi(\al_ix)=\sum_{i=1}^n a_i \phi(\al_i)\phi(x)
\end{equation}
and
\begin{equation}
\label{keykey}
 \tilde{c}\cdot y=\sum_{i=1}^n a_i \phi(\be_i y)=\sum_{i=1}^n a_i \phi(\be_i)\phi(y)
\end{equation}
for any $x, y\in K$. Choosing $x=y=1$:
\begin{equation}
\label{key2}
\sum_{i=1}^n a_i\phi(\al_i)=\sum_{i=1}^n a_i\phi (\be_i)=\tilde{c}.
\end{equation}
If $\tilde{c}=0$ then equations \eqref{key1}, \eqref{keykey} and \eqref{key2} give that $\phi\colon \cc \to \cc$ is a solution of the homogeneous equation
$$\sum_{i=1}^n a_i f(\al_ix+\be_iy)=0\ \ (x, y\in \cc).$$
If $\tilde{c}\neq 0$  then equations \eqref{key1} and \eqref{key2} give that $\phi(x)=x$ ($x\in K$).
The converse statement is clear.  
\hfill $\square$

\begin{cor} If $\tilde{c}\neq 0$ for an exponential function in $S_1^*$ then the space of the additive solutions on $K$ is
$$\frac{c}{\tilde{c}}\cdot x+S_1^{0}.$$
\end{cor}

\begin{cor}
\label{corpa} If there are no automorphisms satisfying
$$\sum_{i=1}^n a_i\phi(\al_i)=\sum_{i=1}^n a_i\phi (\be_i)=0$$
then $S_1^0$ is trivial for any finitely generated field $K\subset \cc$ containing the para\-meters $\al_i$ and $\be_i$ $(i=1, \ldots, n)$ and the only nonzero additive solution of \eqref{e2} on $\cc$ must be the proportional of the identity function: $f(x)=c'\cdot x$ where $c'=c/\tilde{c}$ 
provided that 
$$\sum_{i=1}^n a_i\al_i=\sum_{i=1}^n a_i\be_i=\tilde{c}\neq 0.$$
\end{cor}

\begin{remark}\label{rem1}{\rm 
If $\tilde{c}=0$ for any exponential function in $S_1^*$ then the exponentials give only translation parts  in the solution of the inhomogeneous equation on $K$ and we need to apply spectral synthesis in the variety $S_1^*$ to decide the existence of a nonzero particular solution of the inhomogeneous equation on finitely generated fields containing the parameters $\al_i$ and $\be_i$ $(i=1, \ldots, n)$; see \cite{KV1}.}
\end{remark}

\begin{remark}\label{rem2}
{\rm The size of the subspace $S_1^0$ depends on the existence of 
an automorphism $\phi\colon\mathbb{C}\to \mathbb{C}$ satisfying 
$$\sum_{i=1}^n a_i\phi(\al_i)=0 \ \ \textrm{and} \ \ \sum_{i=1}^n a_i\phi (\be_i)=0.$$
The so-called characteristic polynomial method helps us to investigate such an existence problem in terms of polynomials whose coefficients depend algeb\-raically on the parameters $\al_i$ and $\be_i$ $(i=1, \ldots, n)$; \cite{V1}, see also \cite{KV}, \cite{VV-1}, \cite{V2} and \cite{VV}.}
\end{remark}

\begin{example} \label{exm1} {\rm It is not too hard to find examples such that $S_1^0$ is trivial, see e.g. section 3 in \cite{V1} and section 5 for some explicite examples. In what follows we present only the extremal cases of the characteristic polynomial method: first of all reduce equation $\displaystyle \sum_{i=1}^n a_i\phi(\al_i)=0$ to
\begin{equation}\label{reduced}
1+\sum_{i=1}^{n-1} a_i' \phi(\al_i')=0,\ \ \textrm{where}\ \ a_i'=a_i/a_n\ \ \textrm{and}\ \ \alpha_i'=\alpha_i/\alpha_n,
\end{equation}
$i=1, \ldots, n-1$. 
\begin{itemize}
\item if the parameters $\al_1', \ldots, \al_{n-1}'$ form an algebraically independent system over the rationals then the characteristic polynomial is
$$\mathcal{P}_1(x_1, \ldots, x_{n-1})=1+\sum_{i=1}^{n-1} a_i' x_i.$$
\end{itemize}
The sufficient and necessary condition for the existence of $\phi$ satisfying \eqref{reduced} is that $\mathcal{P}_1$ has a root with algebraically independent coordinates over the rationals. It happens if and only if one of the parameters $a_1'$, $\ldots$, $a_{n-1}'$ is trans\-cendent;  \cite{VV-1}, see also \cite{V2} and \cite{VV}. 
\begin{itemize} 
\item If all the parameters $\al_1', \ldots, \al_{n-1}'$ are algebraic over $\qq$ then the extension of the rationals with $\al_1', \ldots, \al_{n-1}'$ is a simple algebraic extension by an algebraic element $u$ of finite degree. The elements of the extended field can be written as 
$$\al_i'=p_i(u),\ \textrm{where}\ p_i\in \mathbb{Q}[x],$$
$i=1, \ldots, n-1$. This means that the action of $\phi$ is uniquely determined by $\phi(u)$. Consider the finite product
$$\mathcal{F}_1(a_1', \ldots, a_{n-1}')=\prod_{h}\left(1+a_1' p_1(s_h)+\ldots+a_{n-1}'p_{n-1}(s_h)\right),$$
where $s_h$ runs through the finitely many algebraic conjugates of the element $u$. 
\end{itemize}
Using the fundamental theorem of the symmetric polynomials it can be proved that the right hand side is a rational polynomial of the variables $a_1'$, $\ldots$, $a_{n-1}'$: the coefficients are the polynomials of
$\displaystyle{\sum s_h}$, $\displaystyle{\sum s_h\cdot s_k}$ and so on. Its vanishing is the sufficient and necessary condition for the existence of $\phi$ satisfying \eqref{reduced}. The inverse automorphism allows us to change the role of the parameters to unify the polynomial for the family of the parameters $\alpha_i$'s and $\beta_i$'s:
\begin{equation}\label{unify}
\sum_{i=1}^n \phi^{-1}(a_i)\al_i=0 \ \ \textrm{and} \ \ \sum_{i=1}^n \phi^{-1}(a_i)\be_i=0.
\end{equation}}
\end{example}

\section{Applications of spectral analysis II}

In what follows we generalize the previous theorem for the nonzero monomial solutions of higher degree $p>1$ of equation \eqref{ek}

\subsection{The case of higher order monomial solutions}

Suppose that 
$$f(x):=F_p(x, \ldots, x)$$
is a nonzero monomial solution of degree $p>1$ of equation \eqref{ek}, where $F_p$ is a symmetric, $p$-additive function, i.e. $f(e)\neq 0$ for some $e\in \cc$. Let $K=\mathbb{Q} (\al_1,\ldots , \al_n, \be_1, \ldots ,\be_n, e)$ be the extension of $\mathbb{Q}$ by the
complex numbers $\al_i, \be_i$ and $e$ ($i=1,\ldots ,n $). In the sense of Lemma \ref{lsp}, $S_p^*$ is a variety in $\cc^{G^*}$, where 
$$G^*:=\underbrace{K^*\times \ldots \times K^*}_{\textrm{p}-times}.$$
We have that $S_p^*\neq \{0\}$ because $F_p|_{G^*}\in S_p^*$ and $F_p(e, \ldots, e)=f(e)\neq 0$. Since $G^*$ is countable we find, by Theorem \ref{tspanal}, that $S_p^*$ contains an exponential element $\phi$:
$$\phi(x_1y_1, \ldots, x_py_p)=\phi(x_1, \ldots, x_p)\cdot \phi(y_1, \ldots, y_p).$$
Using the decomposition formula
\begin{equation}
\label{decomposition}
\phi(x_1, \ldots, x_p)=\phi(x_1, 1, \ldots,1)\cdot \phi(1, x_2,1 \ldots,1)\cdot \ldots \cdot \phi(1, \ldots, 1, x_p)
\end{equation}
the exponential element can be written as the product
$$\phi=\phi_1 \cdot \ldots \cdot \phi_p,$$
where $\phi_j\colon K^* \to \cc$ ($j=1, \ldots, p$) are exponentials in $\cc^{K^*}$. By Lemma \ref{l30}, each of them can be extended to an automorphism of $\mathbb{C}$. According to the definition of $S_p^*$ (subsection 2.2) the diagonalization
$$\diag \phi (x):=\phi_1(x)\cdot \ldots \cdot \phi_p(x)$$
is the solution of 
\begin{equation}\label{ektilde}
\sum_{i=1}^n a_i \tilde{f}(\al_ix+\be_iy)=\tilde{c}_p\cdot \sum_{l=0}^p x^l y^{p-l} \ \ (x, y\in K)
\end{equation}
for some $\tilde{c}_p\in \cc$. Especially 
\begin{equation}\label{ekkk}
\tilde{c}_p\cdot x^p=\sum_{i=1}^n a_i\phi(\al_ix, \ldots, \al_ix)=\sum_{i=1}^n a_i\phi(\al_i, \ldots, \al_i)\phi(x,\ldots, x)
\end{equation}
and
\begin{equation}\label{ekkkk}
\tilde{c}_p\cdot y^p=\sum_{i=1}^n a_i\phi(\be_i y, \ldots, \be_i y)=\sum_{i=1}^n a_i\phi(\be_i, \ldots, \be_i)\phi(y,\ldots, y).
\end{equation}
Choosing $x=y=1$  
\begin{equation}
\label{ekkkkk}
\sum_{i=1}^n a_i\phi(\al_i, \ldots, \al_i)=\sum_{i=1}^n a_i\phi(\be_i, \ldots, \be_i)=\tilde{c}_p.
\end{equation}
On the other hand 
\begin{equation}
\label{key3}
\begin{split}
&\tilde{c}_p\cdot x^ly^{p-l}=\sum_{i=1}^n a_i \binom{p}{l} \phi(\underbrace{\al_i x, \ldots, \al_i x}_{l \ \textrm{times}}, \be_i y, \ldots, \be_i y)=\\
&\sum_{i=1}^n a_i \binom{p}{l} \phi(\underbrace{\al_i, \ldots, \al_i}_{l \ \textrm{times}}, \be_i, \ldots, \be_i )\phi(\underbrace{x, \ldots, x}_{l \ \textrm{times}}, y, \ldots, y).
\end{split}
\end{equation}
Choosing $x=y=1$
\begin{equation}\label{ekkkkkk}
\sum_{i=1}^n a_i \binom{p}{l} \phi(\underbrace{\al_i, \ldots, \al_i}_{l \ \textrm{times}}, \be_i, \ldots, \be_i )= \tilde{c}_p
\end{equation}
and equation (\ref{ekkkkkk}) holds for any permutation $\sigma$ of the indices, i.e. 
\begin{equation}\label{ekkkkkkk}
\sum_{i=1}^n a_i \binom{p}{l} \phi_{\sigma(1)}(\al_i)\cdot \ldots \cdot \phi_{\sigma(l)}(\al_i)\cdot \phi_{\sigma(l+1)}(\be_i)\cdot \ldots \cdot \phi_{\sigma(p)}(\be_i)= \tilde{c}_p.
\end{equation}
If $\tilde{c}_p=0$ then equations \eqref{ekkk}-\eqref{ekkkkkkk} give that $\diag \phi$ 
is the solution of the homogeneous equation
$$\sum_{i=1}^n a_i f(\al_ix+\be_iy)=0\ \ (x, y\in \cc).$$
If $\tilde{c}_p\neq 0$ then equations \eqref{key3} and \eqref{ekkkkkk} imply that 
$$x^l y^{p-l}=\phi(\underbrace{x, \ldots, x}_{l \ -\ \textrm{times}}, y, \ldots, y).$$
Choosing $l=1$ and $y_2=\ldots y_p=1$ we can conclude that $\phi_1(x)=x$ and so on:
$\phi_1(x)=\ldots=\phi_p(x)=x$ $(x\in K)$. Theorem \ref{t1o1} has the following analogue.  

\begin{theorem}\label{t1ok}
The existence of a nonzero monomial solution of degree $p>1$ of \eqref{ek} implies that there exist a finitely generated subfield $K\subset \cc$ containing $\al_i$ and $\be_i$ $(i=1, \ldots, n)$ and some automorphisms $\phi_i\colon \cc\to \cc$ $(i=1, \ldots, p)$ as the extensions of the functions in the decomposition formula \eqref{decomposition} of an exponential function $\phi$ in $S_p^*$ such that 
\begin{equation}\label{ektildeplus}
\sum_{i=1}^n a_i \diag \phi (\al_ix+\be_iy)=\tilde{c}_p\cdot \sum_{l=0}^p x^l y^{p-l} \qquad (x,y\in K),
\end{equation}
for some $\tilde{c}\in \cc$. Especially
\begin{equation}\label{eAk}
\begin{split}
&\sum_{i=1}^n a_i\phi_1(\al_i)\cdot \ldots \cdot \phi_p(\al_i)=\sum_{i=1}^n a_i\phi_1(\be_i)\cdot \ldots \cdot \phi_p(\be_i)=\\
&\sum_{i=1}^n a_i \binom{p}{l} \phi_{\sigma(1)}(\al_i)\cdot \ldots \cdot \phi_{\sigma(l)}(\al_i)\cdot \phi_{\sigma(l+1)}(\be_i)\cdot \ldots \cdot \phi_{\sigma(p)}(\be_i)= \tilde{c}_p,
\end{split}
\end{equation}
where $l=1, \ldots, p-1$ and $\sigma$ is an arbitrary permutation of the indices. If $\tilde{c}_p=0$ then
\begin{equation}\label{ektilde1}
\sum_{i=1}^n a_i \diag \phi (\al_ix+\be_iy)=0 \qquad (x, y\in \cc),
\end{equation}
i.e. $\diag \phi$ is the solution of the homogeneous equation on $\cc$. If $\tilde{c}_p\neq 0$ then $\ph_1(x)=\ldots=\phi_p(x)=x$ $(x\in K)$ and
\begin{equation}
\label{simplified1}
\sum_{i=1}^n a_i \al_i^p=\sum_{i=1}^n a_i \be_i^p=\sum_{i=1}^n a_i\binom{p}{l} \al_i^l \be_i^{p-l}=\tilde{c}_p\neq 0\ \ (l=1, \ldots, p-1).
\end{equation}
Conversely, if \eqref{simplified1} holds then $f(x):=(c/\tilde{c}_p)\cdot x^p$ is a nonzero particular monomial solution of degree $p$ of \eqref{ek} on $\cc$. 
\end{theorem}

\begin{cor} If $\tilde{c}_p\neq 0$ for an exponential function in $S_p^*$ then the space of the monomial solutions of degree $p$ on $K$ is 
$$\frac{c}{\tilde{c}_p}\cdot x^p+\ \textrm{diag}\ S_p^{0},$$
where $\diag$ means the diagonalizations of the elements of the set.  
\end{cor}

\begin{cor}
\label{corpapa} If there are no automorphisms satisfying 
\begin{equation}
\begin{split}
&\sum_{i=1}^n a_i\phi_{1}(\al_i)\cdot \ldots \cdot \phi_{p}(\al_i)=\sum_{i=1}^n a_i \phi_{1}(\be_i)\cdot \ldots \cdot \phi_{p}(\be_i)=\\
&\sum_{i=1}^n a_i \binom{p}{l} \phi_{\sigma(1)}(\al_i)\cdot \ldots \cdot \phi_{\sigma(l)}(\al_i)\cdot \phi_{\sigma(l+1)}(\be_i)\cdot \ldots \cdot \phi_{\sigma(p)}(\be_i)=0,
\end{split}
\end{equation}
then $S_p^0$ is trivial for any finitely generated field $K\subset \cc$ containing the parameters $\al_i$ and $\be_i$ $(i=1, \ldots, n)$ and the only nonzero monomial solution of degree $p$ of \eqref{ek} on $\cc$  must be the proportional of the $p$th power function:
$$f(x)=c'\cdot x^p,\ \ \textrm{where}\ \ c'=c/\tilde{c}_p$$
provided that 
$$\sum_{i=1}^n a_i \al_i^p=\sum_{i=1}^n a_i \be_i^p=\sum_{i=1}^n a_i\binom{p}{l} \al_i^l \be_i^{p-l}=\tilde{c}_p\neq 0\ \ (l=1, \ldots, p-1).$$
\end{cor}

\begin{remark}\label{rem3}{\rm 
If $\tilde{c}_p=0$ for any exponential function in $S_p^*$ then the diagonalizations of the exponentials give only translation parts in the solution of the inhomogeneous equation on $K$ and we need to apply spectral synthesis in the variety $S_p^*$ to decide the existence of a nonzero particular monomial solution of degree $p$ of the inhomogeneous equation on finitely generated fields containing the parameters $\al_i$ and $\be_i$ $(i=1, \ldots, n)$; see \cite{KV1}.}
\end{remark}

\begin{remark}
{\rm The size of the subspace $S_p^0$ depends on the existence of automorphisms $\phi_1, \ldots, \phi_p\colon\mathbb{C}\to\mathbb{C}$ satisfying 
\begin{equation}
\begin{split}
&\sum_{i=1}^n a_i\phi_{1}(\al_i)\cdot \ldots \cdot \phi_{p}(\al_i)=\sum_{i=1}^n a_i \phi_{1}(\be_i)\cdot \ldots \cdot \phi_{p}(\be_i)=\\
&\sum_{i=1}^n a_i \binom{p}{l} \phi_{\sigma(1)}(\al_i)\cdot \ldots \cdot \phi_{\sigma(l)}(\al_i)\cdot \phi_{\sigma(l+1)}(\be_i)\cdot \ldots \cdot \phi_{\sigma(p)}(\be_i)=0,
\end{split}
\end{equation}
where $l=1, \ldots, p-1$ and $\sigma$ is an arbitrary permutation of the indices. The so-called characteristic polynomial method helps us to investigate such an existence problem in terms of polynomials whose coefficients depend algeb\-raically on the parameters $\al_i$ and $\be_i$ $(i=1, \ldots, n)$; \cite{V1}, see also \cite{KV}, \cite{VV-1}, \cite{V2} and \cite{VV}.}
\end{remark}

\begin{example}
\label{hom}
{\rm It is not too hard to find examples such that $S_p^0$ is trivial, see e.g. a necessary and sufficient condition for $S_2^{0}\neq \{0\}$ in \cite{KVV} provided that $\be_i=1-\al_i$ and $i=1, \ldots, 4$; see also \cite{VV0}. The special choice $\be_i=1-\al_i$ $(i=1, \ldots, n)$ allows us to conclude a descending tendency in the space of the solutions of the homogeneous equation in general: the existence of a nonzero monomial solution of degree $p>1$ implies the existence of nonzero monomial solutions of degree $p-1, \ldots, 1$ too; see \cite{KV} and \cite{KVV}. Therefore $S_1^0=\{0\}$ implies that $S_p^0=\{0\}$, cf. Example \ref{exm1} and section 5.}
\end{example}

\begin{example}\rm{Theorem \ref{t1ok} motivates the problem of the existence of a monomial solution of the form $\phi^p$, where $\phi$ is an automorphism of $\cc$. It is solved in Proposition 3.9 of \cite{KV} by constructing a counterexample: there is a linear functional equation with a solution of degree two but there is no any solution of the form $\phi^2$, where $\phi$ is an automorphism of $\cc$.}
\end{example}

\section{Examples}

The common feature of the following examples is that for any $p=1, \ldots, 2n-1$ the space $S_p^0$ belonging to the restricitions of the solutions of the homogeneous equation is trivial in case of any finitely generated field $K\subset \cc$ containing the parameters $\al_i$ and $\be_i$ $(i=1, \ldots, n)$. This means that each of them contains only the identically zero function and the inhomogeneous equation can be solved by spectral analysis; cf. Corollary \ref{corpa} and Corollary \ref{corpapa}. Therefore we present all solutions of the following functional equations defined on $\rr$ or $\cc$. The examples show how our method is working in explicite cases: the selection is based on our motivating papers \cite{KKSZ08}, \cite{KKSZ} and \cite{KKSZW}. In what follows we suppose that conditions of Lemma \ref{l01} are satisfied. 

\subsection{The first example} In their paper \cite{KKSZ08} the authors consider functional equation
\begin{equation}
\label{szostok}
F(y)-F(x)=(y-x)\left(f(\al x+\be y)+f(\be x+\al y)\right),
\end{equation}
i.e. $n=2,\ a_1=a_2=1,\ \al_1=\al,\ \al_2=\be, \ \be_1=\be,\ \be_2=\al.$ Using Lemma \ref{l01}, Lemma \ref{l1} and Theorem \ref{KKT} the solutions are generalized polynomials of degree at most $3$. Therefore we should solve the following functional equations:
\begin{equation}
\label{KKSZ0}
f(\al x+\be y)+f(\be x+\al y)=c_0,
\end{equation}
\begin{equation}
\label{KKSZ1}
f(\al x+\be y)+f(\be x+\al y)=c_1\cdot (x+y),
\end{equation}
\begin{equation}
\label{KKSZ2}
f(\al x+\be y)+f(\be x+\al y)=c_2\ \cdot (x^2+x\cdot y+y^2),
\end{equation}
\begin{equation}
\label{KKSZ3}
f(\al x+\be y)+f(\be x+\al y)=c_3\ \cdot (x^3+x^2\cdot y+x\cdot y^2+y^3)
\end{equation}
for some constants belonging to the possible values of $p=0, 1, 2, 3$. Equation \eqref{KKSZ0} is trivial; the solutions are the constant functions (monomial terms of degree zero). 

\begin{lemma}
For any $p=1, 2, 3$ we have that $S_p^0$ belonging to the solutions in the homogeneous case $\tilde{c}_p=0$ is trivial for any finitely generated subfield $K\subset \cc$ containing the parameters $\al_i$ and $\be_i$ $(i=1, 2)$.
\end{lemma}

\proof
In case of $p=1$ the non-zero element in $S_1^0$ implies the existence of an automorhism $\phi\colon \cc \to \cc$ such that
$$\sum_{i=1}^2 a_i\phi(\al_i)=\sum_{i=1}^2 a_i\phi(\be_i)=0\ \ \Rightarrow\ \ \phi(\al)+\phi(\be)=0,$$
i.e. $\al=-\be$ which contradicts to condition ($3$) in Lemma \ref{l01}. In case of $p=2$ the non-zero element in $S_2^{0}$ implies the existence of automorphisms $\phi_1,\ \phi_2 \colon \cc \to \cc$ such that
$$\sum_{i=1}^2 a_i\phi_1(\al_i)\phi_2(\al_i)=\sum_{i=1}^2a_i\phi_1(\be_i)\phi_2(\be_i)=\sum_{i=1}^2a_i \binom{2}{1} \phi_{\sigma(1)}(\al_i)\phi_{\sigma(2)}(\be_i)=0$$
for any permutation $\sigma$ of the indices; see Theorem \ref{t1ok}. Especially this system of equations reduces to
$$\phi_1(\al)\phi_2(\al)+\phi_1(\be)\phi_2(\be)=2\phi_1(\al)\phi_2(\beta)+2\phi_1(\beta)\phi_2(\al)=0$$
because of the symmetric roles of $\al$ and $\be$.
Introducing the notations
$$\omega_i=\phi_i\left(\frac{\al}{\be}\right),\ \ \textrm{where}\ \ i=1, 2$$
it follows that
\begin{equation}
\omega_1\cdot \omega_2+1=2(\omega_1+\omega_2)=0
\end{equation}
and, consequently $\omega_1^2=1$, i.e. $\al=\pm \be$ which contradicts to condition ($3$) in Lemma \ref{l01}. In case of $p=3$ the non-zero element in $S_3^0$ implies the existence of automorphisms $\phi_1, \phi_2, \phi_3\colon \cc \to \cc$ such that
$$\sum_{i=1}^2 a_i\phi_1(\al_i)\phi_2(\al_i)\phi_3(\al_i)=\sum_{i=1}^2a_i\phi_1(\be_i)\phi_2(\be_i)\phi_3(\be_i)=$$
$$\sum_{i=1}^2a_i \binom{3}{1} \phi_{\sigma(1)}(\al_i)\phi_{\sigma(2)}(\be_i)\phi_{\sigma (3)}(\be_i)=$$
$$\sum_{i=1}^2a_i \binom{3}{2} \phi_{\sigma(1)}(\al_i)\phi_{\sigma(2)}(\al_i)\phi_{\sigma(3)}(\be_i)=0$$
for any permutation $\sigma$ of the indices. Especially this system of equations reduces to
\begin{equation}
\begin{split}
&\omega_1\cdot \omega_2\cdot \omega_3+1=3\left(\omega_1+\omega_2\cdot \omega_3\right)=3\left(\omega_3+\omega_1\cdot \omega_2\right)=\\
&3\left(\omega_2+\omega_1\cdot \omega_3\right)=0\\
\end{split}
\end{equation}
because of the symmetric roles of $\al$ and $\be$, where 
$$\omega_i=\phi_i\left(\frac{\al}{\be}\right)\ \ \textrm{and}\ \ i=1, 2, 3.$$
Therefore $\omega_1^2=1$, i.e. $\al=\pm \be$ which contradicts to condition ($3$) in Lemma \ref{l01}. 
\hfill $\square$

In the sense of the previous lemma each monomial term in the solution of equation \eqref{szostok} must be the proportional of $x^p$ ($p=1, 2, 3$) provided that the parameters satisfy conditions of Theorem \ref{t1o1} and Theorem \ref{t1ok} under the choice $\phi(x)=x$, $\phi_1(x)=\phi_2(x)=x$ and $\phi_1(x)=\phi_2(x)=\phi_3(x)=x$, respectively:
\begin{itemize} 
\item[(i)] there is an additive (monomial term of degree 1) term in the solution if and only if
$$\sum_{i=1}^2 a_i \al_i=\sum_{i=1}^2 a_i \be_i=\tilde{c}_1 \neq 0\ \ \Rightarrow\ \ \al+\be=\tilde{c}_1 \neq 0,$$
It is obviously true because of condition ($3$) in Lemma \ref{l01}. 
\item[(ii)] There is a monomial term of degree $2$ in the solution if and only if  
$$\sum_{i=1}^2 a_i \al_i^2=\sum_{i=1}^2 a_i \be_i^2=\sum_{i=1}^2 a_i \binom{2}{1} \al_i \be_i=\tilde{c}_2\neq 0 \ \ \Rightarrow\ \ $$
$$\al^2+\be^2=2\left(\al \cdot \be+\be \cdot \al \right)=\tilde{c}_2 \neq 0,$$
i.e. $\al^2+\be^2=4\al \be \neq 0$.
\item[(iii)] There is a monomial term of degree $3$ in the solution if and only if 
$$\sum_{i=1}^2 a_i \al_i^3=\sum_{i=1}^2a_i \be_i^3=\sum_{i=1}^2a_i \binom{3}{1} \al_i\be_i^2=\sum_{i=1}^2a_i \binom{3}{2} \al_i^2 \be_i=\tilde{c}_3\neq 0\ \ \Rightarrow$$ 
$$\al^3+\be^3=3\left(\al^2 \cdot \be+\be^2 \cdot \al \right)=\tilde{c}_3 \neq 0,$$
i.e. $\al^2-\al \be +\be^2=3\al \be\neq 0 \ \ \Rightarrow\ \ \al^2+\be^2=4\al \be \neq 0$.
\end{itemize}
Therefore the spectral analysis presents the results of Theorem 2 (v) and (vi) in \cite{KKSZ08}. Note that condition $\al^2+\be^2=4\al \be$ implies the fraction $\al/\be$ to be $2\pm \sqrt{3}$. The discussion of the exceptional cases $\al=\pm \be$ or $\be=0$ can be found in Theorem 2 (i) - (iv) of \cite{KKSZ08}. In these cases condition ($3$) of Lemma \ref{l01} is taking to fail.  

\begin{remark}\label{fail} \rm{If condition ($3$) in Lemma 1.1 is taking to fail then there can be more general solutions: under the choice $\al=-\be$, equation \eqref{szostok} is satisfied by any odd functions with $F=\textrm{\ constant}$; cf. Remark 3 in \cite{KKSZ}}.
\end{remark}

\subsection{The second example} In their paper \cite{KKSZ} the authors consider functional equation
\begin{equation}
\label{szostok1}
F(y)-F(x)=(y-x)\sum_{i=1}^n a_i f(\al_i x+(1-\al_i)y),
\end{equation}
where $\sumina \neq 0$ and $\be_i=1-\al_i$ ($i=1, \ldots, n$). In what follows we also suppose that the parameters $\al_i$'s are different to satisfy condition ($3$) in Lemma \ref{l01}; see remark \ref{fail}.
Using Lemma \ref{l01}, Lemma \ref{l1} and Theorem \ref{KKT} the solutions are generalized polynomials of degree at most $2n-1$. Therefore we should solve functional equations of type 
\begin{equation}
\label{KKSZp}
\sum_{i=1}^n a_i f(\al_i x+(1-\al_i)y)=c_p\sum_{l=0}^p x^p y^{p-l},
\end{equation}
where $p=0, 1, \ldots, 2n-1$ and $\sum_{i}^n a_i\neq 0$. The following Lemma shows that condition $\sum_{i}^n a_i\neq 0$ means the disqualification of non-zero solutions of the homogeneous equation.

\begin{lemma} \label{hom1}
For any $p\geq 1$ we have that $S_p^0$ belonging to the solutions in the homogeneous case $\tilde{c}_p=0$ is trivial  for any finitely generated subfield $K\subset \cc$ containing the parameters $\al_i$ and $\be_i$ $(i=1, \ldots,n)$. 
\end{lemma} 

\proof
Substituting $x=y$ in equation  
\begin{equation}
\label{KKSZp0}
\sum_{i=1}^n a_i f(\al_i x+(1-\al_i)y)=0
\end{equation}
we have that $f(x)=0$ because of $\sum_{i=1}^na_i\neq 0$.
\hfill $\square$

\begin{remark} \rm{Equation \eqref{KKSZp0} is widely studied under the condition $\sum_{i=1}^n a_i=0$; \cite{KV}, \cite{KVV}, see also \cite{VV0}. The special choice $\be_i=1-\al_i$ $(i=1, \ldots, n)$ allows us to conclude a descending tendency in the space of the solutions of the homogeneous equation in general: the existence of a nonzero monomial solution of degree $p>1$ implies the existence of nonzero monomial solutions of degree $p-1, \ldots, 1$ too; see \cite{KV} and \cite{KVV}, Example \ref{hom}. This means that condition $\sumina \neq 0$ is essential to conclude Lemma \ref{hom1}; cf. Remark 2 in \cite{KKSZ}.}  
\end{remark}

\begin{remark}\rm{The proof of Lemma \ref{hom1} is also working in a more general situation according to the rational homeogenity of $p$-additive functions: let us choose rational numbers $r_1, \ldots r_n$
such that the elements
$${\bf r}_1=(r_1, \ldots, r_n), \ldots, {\bf{r}}_{2n-1}:=(r_1^{2n-1}, \ldots,r_n^{2n-1})$$
belong to the half-space defined by $\langle {\mathbf{\alpha}}, {\bf{x}} \rangle >0$, where ${\mathbf{\alpha}}=(\al_1, \ldots, \al_n).$ If $\be_i=r_i-\al_i$ ($i=1, \ldots, n$) then, by repeating the steps of the proof of Lemma \ref{hom1}, it can be easily seen that $S_p^{0}$ is trivial for any $p=1, \ldots, 2n-1$. The construction is the general version of Remark 4 in \cite{KKSZ}.}
\end{remark}

In the sense of the previous lemma each monomial term in the solution of equation \eqref{szostok1} must be the proportional of $x^p$ ($p=1, \ldots 2n-1)$ (see also Theorem 1 in \cite{KKSZ}). Now we apply conditions of Theorem \ref{t1o1} and Theorem \ref{t1ok} under the choice $\phi(x)=x$, $\phi_1(x)=\phi_2(x)=x$, $\ldots$, $\phi_1(x)=\phi_2(x)=\ldots=\phi_p(x)=x$, respectively:
\begin{equation}
\label{descending}
\sum_{i=1}^n a_i \al_i^p=\sumina (1-\al_i)^p=\sumina \binom{p}{l}\al_i^{l}(1-\al_i)^{p-l}=\tilde{c}_p\neq 0,
\end{equation}
where $l=1, \ldots, p-1$. Under the choice $l=p-1$ we have that
$$\sumina \binom{p}{p-1}\al_i^{p-1}(1-\al_i)=\tilde{c}_p\ \ \Rightarrow \ \ \sumina \al_i^{p-1}=\frac{p+1}{p}\tilde{c}_p=:\tilde{c}_{p-1}.$$
If $l=p-2$ then
$$\sumina \binom{p}{p-2}\al_i^{p-2}(1-\al_i)^{2}=\tilde{c}_p\ \ \Rightarrow\ \ \sumina \al_i^{p-2}=\frac{p+1}{p-1}\tilde{c}_p.$$
Under the choice $l=p-3$:
$$\sumina \binom{p}{p-3}\al_i^{p-3}(1-\al_i)^{3}=\tilde{c}_p\ \ \Rightarrow \ \ \sumina \al_i^{p-3}=\frac{p+1}{p-2}\tilde{c}_p$$
and so on. Therefore system \eqref{descending} is equivalent to
\begin{equation}
\label{descending1}
\begin{split}
&\sum_{i=1}^n a_i \al_i^p=\tilde{c}_p\\
&\sum_{i=1}^n a_i \al_i^{p-1}=\frac{p+1}{p}\tilde{c}_{p}=\tilde{c}_{p-1}\\
&\sum_{i=1}^n a_i \al_i^{p-2}=\frac{p+1}{p-1}\tilde{c}_p=\frac{p}{p-1}\tilde{c}_{p-1}\\
&\sum_{i=1}^n a_i \al_i^{p-3}=\frac{p+1}{p-2}\tilde{c}_p=\frac{p}{p-2}\tilde{c}_{p-1}\\
& ...\\
\end{split}
\end{equation}
System (\ref{descending1}) shows a kind of descending tendency in this case. If the solution of \eqref{szostok1} contains a monomial term of degree $p$ with $\tilde{c}_p \neq 0$ then we also have a monomial term of degree $p-1$ with $\tilde{c}_{p-1}\neq 0$.  

\subsection{The third example} Using some linear substitutions, it can be easily seen that the homogeneous equation \eqref{KKSZp0} is equivalent to 
\begin{equation}
\label{subslinear}
\sum_{i=1}^n a_if(x+\be_i y)=0
\end{equation}
for some parameters $\be_i$'s. Therefore we are motivated to investigate the inhomogeneous case under the choice of the parameters $\al_i=1$ ($i=1, \ldots, n$). Consider functional equation
\begin{equation}
\label{gergo1}
F(y)-F(x)=(y-x)\sum_{i=1}^n a_i f(x+\beta_i y),
\end{equation}
where $\sumina \neq 0$ and $\beta_i$'s are different parameters. According to the equivalence of \eqref{subslinear} and \eqref{KKSZp0} in the homogeneous case there is no need to repeat the investigation of $S_1^{0}, \ldots, S_p^0$. We can directly conclude that each monomial term in the solution of equation \eqref{gergo1} must be the proportional of $x^p$ ($p=1, \ldots 2n-1)$ as the solution of the functional equations of type
\begin{equation}
\label{subsKKSZp}
\sum_{i=1}^n a_i f(x+\be_i y)=c_p\sum_{l=0}^p x^p y^{p-l},
\end{equation}
where $p=0, 1, \ldots, 2n-1$ and $\sum_{i}^n a_i\neq 0$. Using conditions of Theorem \ref{t1o1} and Theorem \ref{t1ok} under the choice $\phi(x)=x$, $\phi_1(x)=\phi_2(x)=x$, $\ldots$, $\phi_1(x)=\phi_2(x)=\ldots=\phi_p(x)=x$, respectively:
$$\sum_{i=1}^n a_i=\sumina \be_i^p=\sumina \binom{p}{l}\be_i^{p-l}=\tilde{c}_p\neq 0,$$
where $l=1, \ldots, p-1$. Observe that for any $p=1, \ldots, 2n-1$ we have a universal constant $\tilde{c}_p$ because of $\tilde{c}_p=\sumina$ ($p=1, \ldots, 2n-1$). This means that the solution must be the proportional of $x^p$ for a \emph{uniquely determined value} of the power $p$. Therefore there is no descending tendency and the family of the parameters $\al_i$ and $1-\al_i$ ($i=1, \ldots,n$) or $\al_i=1$ and $\be_i$ ($i=1, \ldots,n$) represent essentially different classes of functional equations in the inhomogeneous case.  

\subsection{The fourth example} In their paper \cite{KKSZW} the authors consider functional equation
\begin{equation}
\label{KKSZW1}
F(y)-F(x)=(y-x)\left(af(x)+(1-a)f(\lambda x+(1-\lambda)y)\right),
\end{equation}
i.e. $n=2$, $a_1=a$, $a_2=1-a$, $\al_1=1$, $\al_2=\lambda$, $\be_1=0$, $\be_2=(1-\lambda)$, where $\lambda\neq 0, 1$ and $a\neq 0$. Using Lemma \ref{l01}, Lemma \ref{l1} and Theorem \ref{KKT} the solutions are generalized polynomials of degree at most $3$. Therefore we should solve functional equations 
\begin{equation}
\label{KKSZW0}
af(x)+(1-a)f(\lambda x+(1-\lambda)y)=c_0,
\end{equation}
\begin{equation}
\label{KKSZW01}
af(x)+(1-a)f(\lambda x+(1-\lambda)y)=c_1\cdot (x+y),
\end{equation}
\begin{equation}
\label{KKSZW02}
af(x)+(1-a)f(\lambda x+(1-\lambda)y)=c_2\ \cdot (x^2+x\cdot y+y^2),
\end{equation}
\begin{equation}
\label{KKSZW03}
af(x)+(1-a)f(\lambda x+(1-\lambda)y)=c_3\ \cdot (x^3+x^2\cdot y+x\cdot y^2+y^3)
\end{equation}
for some constants belonging to the possible values of $p=0, 1, 2, 3$. Equation \eqref{KKSZW0} is trivial; the solutions are the constant functions (monomial terms of degree zero). The translation parts of a non-zero particular solution of the inhomogeneous equations depends on the size of $S_1^0$, $S_2^0$ and $S_3^0$.  We are going to prove that all of them are trivial:
$$\sum_{i=1}^2 a_i \phi(\al_i)=\sum_{i=1}^2 a_i \phi(\be_i)=0$$
implies that 
$$a+(1-a)\phi(\lambda)=(1-a)\left(1-\phi(\lambda)\right)=0,$$
i.e. $\lambda=1$ (which is a contradiction) or $a=1$ and $a=0$ at the same time. The contradictions show that $S_1^0$ contains only the identically zero function. Using that $\phi(1)=1$ and $\phi(0)=0$ for any automorphism of $\cc$ the second order conditions 
$$\sum_{i=1}^2 a_i\phi_1(\al_i)\phi_2(\al_i)=\sum_{i=1}^2a_i\phi_1(\be_i)\phi_2(\be_i)=\sum_{i=1}^2a_i \binom{2}{1} \phi_{\sigma(1)}(\al_i)\phi_{\sigma(2)}(\be_i)=0$$
give the same contradictions as above because of $\lambda\neq 0, 1$. This means that $S_2^0$ is trivial. The third order conditions 
$$\sum_{i=1}^2 a_i\phi_1(\al_i)\phi_2(\al_i)\phi_3(\al_i)=\sum_{i=1}^2a_i\phi_1(\be_i)\phi_2(\be_i)\phi_3(\be_i)=$$
$$\sum_{i=1}^2a_i \binom{3}{1} \phi_{\sigma(1)}(\al_i)\phi_{\sigma(2)}(\be_i)\phi_{\sigma_3}(\be_i)=$$
$$\sum_{i=1}^2a_i \binom{3}{2} \phi_{\sigma(1)}(\al_i)\phi_{\sigma(2)}(\al_i)\phi_{\sigma(3)}(\be_i)=0$$
can be also given as a simple system of equations because of $\be_1=0$. An easy computation shows that $S_3^0$ is trivial. As the next step consider the particular nonzero solutions of equations \eqref{KKSZW01}, \eqref{KKSZW02} and \eqref{KKSZW03}. 
\begin{itemize}
\item[(i)] For the existence of the additive (monomial term of degree $1$) term in the solution it is sufficient and necessary that
$$\sum_{i=1}^2 a_i\al_i=\sum_{i=1}^2 a_i \be_i=\tilde{c}_1\neq 0.$$
Especially $a+(1-a)\lambda=(1-a)(1-\lambda)=\tilde{c}_1 \neq 0,$ i.e. $a\neq 1$ and 
\begin{equation}
\label{descending2}
2\lambda=1-\frac{a}{1-a}
\end{equation}
implies that $a\neq 1/2$ because of $\lambda\neq 0, 1$. 
\item[(ii)] For the existence of the monomial term of degree $2$ in the solution it is sufficient and necessary that
$$\sum_{i=1}^2 a_i \al_i^2=\sum_{i=1}^2 a_i \be_i^2=\sum_{i=1}^2 a_i \binom{2}{1} \al_i \be_i=\tilde{c}_2\neq 0.$$
Especially $a+(1-a)\lambda^2=(1-a)(1-\lambda)^2=2\lambda(1-\lambda)=\tilde{c}_2\neq 0,$
i.e. $\lambda=1/4$ and $a=1/3$. Note that \eqref{descending2} is satisfied under these special values, i.e. we can see the descending tendency too. 
\item[(iii)]  For the existence of the monomial term of degree $3$ in the solution it is sufficient and necessary that
$$\sum_{i=1}^2 a_i \al_i^3=\sum_{i=1}^2a_i \be_i^3=\sum_{i=1}^2a_i \binom{3}{1} \al_i\be_i^2=\sum_{i=1}^2a_i \binom{3}{2} \al_i^2 \be_i=\tilde{c}_3\neq 0.$$ 
Especially
$$a+(1-a)\lambda^3=(1-a)(1-\lambda)^3=3(1-a)\lambda(1-\lambda)^2=$$
$$3(1-a)\lambda^2(1-\lambda)=\tilde{c}_3\neq 0.$$
Therefore
$$\frac{a}{1-a}+\lambda^3=(1-\lambda)^3=3\lambda(1-\lambda)^2=3\lambda^2(1-\lambda).$$
From equation
$$3\lambda(1-\lambda)^2=3\lambda^2(1-\lambda)$$
it follows that $\lambda=1/2$ which contradicts to equation
$$(1-\lambda)^3=3\lambda(1-\lambda)^2$$
and there is no monomial term of degree $3$ in the solution; cf. Theorem 8 in \cite{KKSZW}.
\end{itemize}

\subsection{The fifth example} In their paper \cite{KKSZW} the autors consider functional equation
\begin{equation}
\label{KKSZW2}
F(y)-F(x)=(y-x)\left(af(\lambda x+(1-\lambda)y)+(1-a)f((1-\lambda)x+\lambda y)\right),
\end{equation}
i.e. $n=2$, $a_1=a$, $a_2=1-a$, $\al_1=\lambda$, $\al_2=1-\lambda$, $\be_1=(1-\lambda)$, $\be_2=\lambda$, where $\lambda\neq 0, 1$ and $a\neq 0$. Using Lemma \ref{l01}, Lemma \ref{l1} and Theorem \ref{KKT} the solutions are generalized polynomials of degree at most $3$. Since it is a special case of the second  example with $\sum_{i=1}^2 a_i=1\neq 0$ we can directly conclude that $S_1^0$, $S_2^0$ and $S_3^0$ are trivial. To avoid some unnecessary repetitions we summarize the basic results (see example 5.2):
\begin{itemize}
\item[(i)] for the existence of the additive (monomial term of degree $1$) term in the solution it is sufficient and necessary that
$$a\lambda+(1-a)(1-\lambda)=a(1-\lambda)+(1-a)\lambda=\tilde{c}_1 \neq 0.$$
Equation $a\lambda+(1-a)(1-\lambda)=a(1-\lambda)+(1-a)\lambda$ implies that
$$(1-2a)(1-2\lambda)=0,$$
i.e. $a=1/2$ or $\lambda=1/2$.
\item[(ii)] For the existence of the monomial term of degree $2$ in the solution it is sufficient and necessary that
$$a\lambda^2+(1-a)(1-\lambda)^2=a(1-\lambda)^2+(1-a)\lambda^2=$$
$$2a\lambda(1-\lambda)+2(1-a)(1-\lambda)\lambda=\tilde{c}_2\neq 0.$$
Since we know that the descending tendency holds in the space of the solutions it is enough to check the values $a=1/2$ or $\lambda=1/2$ as we have seen in (i). If $a=1/2$ then we have that
$$6\lambda^2-6\lambda+1=0,\ \ \textrm{i.e.} \ \ \lambda_{12}=\frac{3\pm \sqrt{3}}{6}.$$
In case of $\lambda=1/2$ we have a contradiction:
$$\frac{a}{4}+\frac{1-a}{4}=\frac{a}{4}+\frac{1-a}{4}=\frac{a}{2}+\frac{1-a}{2}=\tilde{c}_2\neq 0.$$
\item[(iii)]  For the existence of the monomial term of degree $3$ in the solution it is sufficient and necessary that
$$a\lambda^3+(1-a)(1-\lambda)^3=a(1-\lambda)^3+(1-a)\lambda^3=$$
$$3a\lambda(1-\lambda)^2+3(1-a)(1-\lambda)\lambda^2=$$
$$3a\lambda^2(1-\lambda)+3(1-a)(1-\lambda)^2\lambda=\tilde{c}_3\neq 0.$$
According to the descending tendency it is enough to check the values
$$a=\frac{1}{2}\ \ \textrm{and}\ \ \lambda_{12}=\frac{3\pm \sqrt{3}}{6}.$$
A direct computation shows that they satisfy the necessary and sufficient conditions; cf. Theorem 9 in \cite{KKSZW}.
\end{itemize}

\section{Concluding remarks} The application of the spectral analysis in the solution of linear functional equations provides a unified method of treatment of the equations in both the homogeneous and inhomogeneous cases. In case of inhomogeneous problems this gives a necessary condition for the existence of the nonzero solutions. The sufficiency also follows in some special cases and the spectral analysis allows us to describe all solutions provided that $S_1^0$, $\ldots$, $S_{2n-1}^0$ are trivial, i.e. they contain only the identically zero function. It depends only on the algebraic properties of the parameters. These algebraic properties are accumulated in the characteristic polynomial of the functional equation; see Remark \ref{rem2} and Example \ref{exm1} (theoretical examples). An explicite example for non-trivial $S_1^0$ can be found in section 3 of  \cite{V1}. In this case we need the so-called spectral synthesis to decide the existence of the nonzero particular solution of the functional equation. On the other hand the spectral synthesis helps us to describe the entire space of the solutions on a large class of finitely generated fields; see \cite{KV1}.

\section{Acknowledgement}
The authors would like to thank to Adrienn Varga for paying our attention to Szostok's problem (the 51st International Symposium on Functional Equations, Rzesz\'{o}w, Poland, June 16-23, 2013).

\end{document}